\documentclass[12pt,russian]{article}
\usepackage{amsfonts,amsmath,amsxtra}
\newtheorem{theorem}{Theorem}
\newtheorem{lemma}{Lemma}

\newenvironment{definition}
{\smallskip\noindent{\bf Definition\/}:}{\smallskip\par}
\newenvironment{proof}
{\noindent{\bf Proof\/}.}{{ $\Box$}\smallskip\par}

\title{On divisorial filtrations on sheaves.}
\author{E.Gorsky.}
\date{}
\begin{document}

\def\eps{\varepsilon}

\maketitle

Moscow State University, Department of Mathematics and Mechanics\footnote{e. mail:
gorsky@mccme.ru}

\sloppy
\begin{abstract}
A notion of Poincar\'e series was introduced in \cite{cdk}. It was developed in \cite{dg} for a multi-index filtration
corresponding to the sequence of blow-ups.
The present
paper suggests the way to generalize the notion of Poincar\'e series to the case
of arbitrary locally free sheaf on the modification of complex plane
$\mathbb{C}^2$.
 This series is expressed through the topological invariants of the sheaf.
For the sheaf of holomorphic
functions the answer coincides with the Poincar\'e series from
\cite{dg}.
\end{abstract}

\section{Introduction}

In \cite{dg} F.Delgado and S.M.Gusein-Zade have computed the
Poincar\'e series of the multi-index filtration defined by a finite
collection of divisorial valuations on the ring
$\mathcal{O}_{\mathbb{C}^2,0}$ of germs of functions of two variables.
Similar to functions, the pull back map lifts holomorphic 1-forms to the
space of modification.

Therefore one could define a
filtration on the space of germs of holomorphic 1-forms on $\mathbb{C}^2.$
This filtration naturally corresponds to a filtration on the space of global
sections of the sheaf of 1-forms on the plane's modification.

Calculating the Poincar\'e series of this filtration seems to be much more
difficult than for the functions.
Hence it is suggested to substitute the space of global sections by the
corresponding sheaf and to calculate Euler characteristics of the
quotient sheaves, organizing them into a generating series.

It is shown below that  the answer for
$\mathcal{O}_{\mathbb{C}^2,0}$ coincides with one from \cite{dg}.
Theorem 1 gives the formula  for an
arbitrary locally free sheaf on the space of modification in terms of
Chern classes of its restriction onto the exceptional lines.

As an example, the series for the sheaf
of 1-forms is calculated in Theorem 2.

By $\underline{v}$ denote the element $(v_1,\ldots,v_s)$ of the lattice
$\mathbb{Z}^s$.
There is a natural partial ordering on $\mathbb{Z}^s$:
$\underline{v}\le\underline{w}$, if every coordinate of $\underline{v}$ is
less or equal to the corresponding coordinate of $\underline{w}$.
For the pair $\underline{v}, \underline{w}\in\mathbb{Z}^s$ let the upper bound
$\sup\{\underline{v}, \underline{w}\}$
be the smallest (according to this ordering) element of $\mathbb{Z}^s$, which is more or equal
 $\underline{v}$ and $\underline{w}$.

\begin{definition} A decreasing
$s$-index filtration on the vector space
$L$ is the family of subspaces
$\{L(\underline{v})|\underline{v}\in
Z^s\}$ such that the following conditions hold:

1) if $\underline{v}_1\le \underline{v}_2$, then
$L(\underline{v}_1)\supset L(\underline{v}_2)$;

2) $L(\underline{v})\cap
L(\underline{w})=L(\sup\{\underline{v},\underline{w}\});$

3) $L(\underline{0})=L\quad (\underline{0}=(0,\ldots,0)).$
\end{definition}

Let $L(\underline{v})$ be a $s$-index filtration on the space
$L$, and all quotient spaces
$L(\underline{v})/L(\underline{v}+\underline{1})$
$(\underline{1}=(1,\ldots,1))$ are finite dimensional. Denote $d(\underline{v})=\dim
L(\underline{v})/L(\underline{v}+\underline{1})$. From (2) and (3) it
follows that for all
$v'_{i_0}<v''_{i_0}\le0$ one has
$$L(v_1,\ldots,v'_{i_0},\ldots,v_s)=L(v_1,\ldots,v''_{i_0},\ldots,v_s),$$
hence the filtration is defined by the set of subspaces
$L(\underline{v})$ with $\underline{v}$ such that all its components are
nonnegative.

Let
$\mathcal{L}=\mathbb{Z}[[t_1,\ldots,t_s,{t_1}^{-1},\ldots,{t_s}^{-1}]]$
be the space of the formal Laurent series of $s$ variables. Elements $\mathcal{L}$
are of the form $\sum_{\underline{v}\in \mathbb{Z}^s}k(\underline{v})\cdot
\underline{t}^{\underline{v}},$ generally speaking, infinite in all
directions. $\mathcal{L}$ is not a ring, but a module over the space of polynomials.
Let
$$Q(t_1,\ldots,t_s)=\sum_{\underline{v}\in \mathbb{Z}^s}d(\underline{v})\cdot
\underline{t}^{\underline{v}}.$$
Since for $v'_{i_0}<v''_{i_0}\le0$ one has

$$d(v_1,\ldots,v'_{i_0},\ldots,v_s)=d(v_1,\ldots,v''_{i_0},\ldots,v_s),$$
the expression
$$P'(t_1,\ldots,t_s)=Q(t_1,\ldots,t_s)\cdot\prod_{i=1}^s(t_i-1)$$ is a power
series, i. e. an element of the subset
$\mathbb{Z}[[t_1,\ldots,t_r]]\subset\mathcal{L}$.

\begin{definition}
We call the series
$$P_L(t_1,\ldots,t_s)={P'(t_1,\ldots,t_r)\over{t_1\cdots t_s -1}}$$
the Poincar\'e series of the multi-index filtration $\{L(\underline{v})\}$
on the space $L$.
\end{definition}
This definition of Poincar\'e series was introduced in
\cite{cdk}.

Let $\pi:(\mathcal{X},\mathcal{D})\rightarrow(\mathbb{C}^2,0)$ be a proper
analytic map which is an isomorphism outside of the origin in $\mathbb{C}^2$
such that $\pi$ is obtained by a sequence of $s$ point blow-ups. Therefore
the exceptional divisor $\mathcal{D}$ is the union of $s$ irreducible components
$E_i$, each of them is isomorphic to the complex projective line.
\begin{lemma}
The pull back map $\pi^{*}$ is an isomorphism between
$H^{0}(\mathbb{C}^2,\Omega^{k}_{\mathbb{C}^2})$ and
$H^{0}(\mathcal{X},\Omega^{k}_{\mathcal{X}})$.
\end{lemma}
\begin{proof}
Let $w$ be a holomorphic $k$-form on $\mathbb{C}^2$ such that $\pi^{*}w=0$.
Then $w=0$ outside the origin, so  $w=0$, and
 $\ker\pi^{*}=0$.

Let $\widetilde{w}$ be a holomorphic $k$-form on the variety $\mathcal{X}$. Since $\pi^{*}$
is an isomorphism outside the exceptional divisor,
there is well-defined holomorphic $k$-form $w=(\pi^{*})^{-1}\widetilde{w}$ on  $\mathbb{C}^2\setminus\{0\}$.
By the Hartogs' theorem it could be continued to an holomorphic $k$-form $w$ on $\mathbb{C}^2$.
Then
$\pi^{*}w$ coincides with $\widetilde{w}$ outside $\mathcal{D}$, so they
coincide in every point of
$\mathcal{X}$. Therefore
$Im\,\pi^{*}=H^{0}(\mathcal{X},\Omega^{k}_{\mathcal{X}}),$ so $\pi^{*}$ is
an isomorphism.
\end{proof}
Let $\mathcal{J}_i$ be the sheaf of ideals of a subscheme $E_i$ in $\mathcal{X}$ and $$d_{i}^{j}=-(E_i\circ
E_j),\,M=(m_{i}^{j})=D^{-1},$$
where $(\circ)$ denotes an intersection number.
Let $\widetilde{E}_i$ be the smooth part of line $E_i$, i. e., the component
$E_i$ without intersection points with all other components of the
exceptional divisor.

Let $\mathcal{G}$ be an arbitrary locally free sheaf on
$\mathcal{X}$.
For every collection $\underline{k}=(k_1,\ldots,k_s)$ of nonnegative numbers
multiplication induces a natural embedding
$\mathcal{G}\otimes\prod_{i=1}^s\mathcal{J}_i^{k_i}\hookrightarrow\mathcal{G}.$
Therefore a multi-index filtration with
subsheaves is defined on $\mathcal{G}$:
$$\mathcal{G}(k_1,\ldots,k_s)=\mathcal{G}\otimes\prod_{i=1}^s\mathcal{J}_i^{k_i}.$$
This filtration defines a multi-index filtration
$H^{0}(\mathcal{X},\mathcal{G}(k_1,\ldots,k_s))$
on the space of the global sections $H^{0}(\mathcal{X},\mathcal{G})$.
Poincar\'e series of the last filtration for
 $\mathcal{G}=\mathcal{O}_{\mathcal{X}}$ was computed in [1]:
$$P_{\mathcal{O}}=\prod_{i}(1-t^{\underline{m}^i})^{-\chi(\widetilde{E}_i)}.$$

Following lemma 1, $H^{0}({\mathcal{X}},\mathcal{O}_{\mathcal{X}}(k_1,\ldots,k_s))$
could be considered as the space of functions $f$ on $\mathbb{C}^2$ such
that $\pi^{*}f$ vanishes on component $E_i$ of the exceptional divisor with
multiplicity higher or equal to $k_i$.
Consider a filtration on the space of 1-forms over $\mathbb{C}^2$ where
subspace with label $\underline{k}$ consists of forms $w$ such that $\pi^{*}w$ vanishes on $E_i$
with multiplicity higher or equal $k_i$.
From lemma 1 it follows that Poincar\'e series of this filtration coincides
with Poincar\'e series of the filtration which is defined above in case
$\mathcal{G}=\Omega^{1}_{\mathcal{X}}$.

Suppose $\pi$ is the unique blow-up.
Let $\omega=f(x,y)dx+g(x,y)dy$ be an 1-form, $k$ is the smallest of orders of functions
 $f$ and $g$ at the origin. Decompose $\omega$ in the sum
$$\omega=\omega_k+\omega_{k+1}+\ldots$$ where $$\omega_m={f_m}dx+{g_m}dy,$$ $f_m$
and
$g_m$ are the homogeneous components of $f$ and $g$ with degree $m$.
Let $\theta$ be the affine coordinate in chart $\{x\neq 0\}$ of the exceptional divisor.
Then $y=\theta\cdot x$, so
\begin{flushleft}
$\pi^{*}\omega=(f_{k}(x,y)+\theta\cdot{g_k}(x,y)+f_{k+1}(x,y)+\theta{g_{k+1}}(x,y))dx+xg_{k}(x,y)dy+...=
x^k(f_{k}(1,\theta)+\theta\cdot{g_k}(1,\theta))dx+x^{k+1}(f_{k+1}(1,\theta)+\theta{g_{k+1}}(1,\theta))dx+g_{k}(1,\theta)dy)+...$
\end{flushleft}
Hence $\pi^{*}\omega$ has multiplicity on $\mathcal{D}$ bigger than $k$ if and only if
 $$f_{k}(x,y)+\theta\cdot{g_k}(x,y)=0,$$ i. e.,
$xf_k+yg_k=0$, so $xf_k=-yg_k.$
Therefore $$f_k=-y\varphi,g_k=x\varphi$$ for a homogeneous polynomial $\varphi(x,y)$ with degree
$k-1$, thus
$$\omega_k=\varphi(xdy-ydx).$$
Then
$$\pi^{*}\omega=x^{k+1}(f_{k+1}(1,\theta)+\theta{g_{k+1}}(1,\theta))dx+g_{k}(1,\theta)dy)+...,$$
and since $g_k\neq 0$, the form $\pi^{*}\omega$ has a multiplicity equal to $k+1$
on
$\mathcal{D}$.

Finally, $F(v)$ consists of forms $$\varphi(xdy-ydx)+fdx+gdy,$$
where $\varphi$ is a homogeneous polynomial with degree $v-2$ and $f$ and $g$
have orders, bigger or equal to $v$. Let us compute $P(t)$. Denote by $d_k=k+1$ the
dimension of the space $D_k$ of homogeneous polynomials of two variables with degree
$k$, $R_k$ be the space of 1-forms  $\varphi(xdy-ydx)$ where $\varphi$ is a
homogeneous polynomial with degree $k-1$. Then $\dim R_k=d_{k-1}$,
$$F(k)/F(k+1)\cong[(D_{k}\oplus D_{k})/R_{k-1}]\oplus R_{k-2},$$
hence $\dim(F(k)/F(k+1))=2d_{k}-d_{k-1}+d_{k-2},$
$$
P(t)={{2-t+t^2}\over(1-t)^2}={1+t\over(1-t)^2}+1.
$$


\section{Geometrical Poincar\'e series }

The space of global sections and its dimension are quite sophisticated
invariants of a sheaf. It is simpler to calculate its Euler characteristic.
Let
$$h(\underline{v})=\chi\biggl(\mathcal{X},\mathcal{G}(\underline{v})/\mathcal{G}(\underline{v}+\underline{1})\biggr).$$
Since the sheaf
$\mathcal{G}(\underline{v})/\mathcal{G}(\underline{v}+\underline{1})$
is supported on $\mathcal{D}$, one could see that $$h(\underline{v})=\chi\biggl(\mathcal{D},\mathcal{G}(\underline{v})/\mathcal{G}(\underline{v}+\underline{1})\biggr).$$
Consider a formal Laurent series
$\widetilde{Q}(t_1,\ldots,t_s)=\sum_{\underline{v}\in \mathbb{Z}^s}h(\underline{v})\cdot\underline{t}^{\underline{v}}$
and denote
$$\widetilde{P}'(t_1,\ldots,t_s)=\widetilde{Q}(t_1,\ldots,t_s)\cdot\prod_{i=1}^s(t_i-1).$$
As above, $\widetilde{P}'$ is a power series, and the series
$$\widetilde{P}_{\mathcal{G}}(t_1,\ldots,t_s)={\widetilde{P}'(t_1,\ldots,t_s)\over{t_1\cdots
t_s-1}}$$ is well-defined.

For $I\subset I_0=\{1,2,\ldots,s\}$ by $|I|$ denote the number of elements in
 $I$, and let $\underline{1}_I$ be the element of $\mathbb{Z}^r$ such that
its
components labeled by numbers from $I$ are equal to 1, and other components are equal to zero.
Denote $$k(\underline{v})=-\sum_{I\subset
I_0}(-1)^{|I|}\chi\biggl(\mathcal{D},\mathcal{G}(\underline{v})/\mathcal{G}(\underline{v}+\underline{1}_I)\biggr).$$
\begin{lemma}
$$\widetilde{P}_\mathcal{G}(\underline{t})=\sum_{\underline{v}\in \mathbb{Z}^s_{\ge0}}k(\underline{v})\cdot
\underline{t}^{\underline{v}}.$$
\end{lemma}
\begin{proof}
The coefficient at the monomial
${\underline{t}}^{\underline{v}}$ in the series
$$\biggl(\sum k(\underline{v})\cdot
\underline{t}^{\underline{v}}\biggr)(t_1\cdots t_s -1)$$ is equal to $$-\sum_{I\subset
I_0}(-1)^{|I|}\chi\biggl(\mathcal{D},\mathcal{G}(\underline{v}-\underline{1})/\mathcal{G}(\underline{v}-\underline{1}+\underline{1}_I)\biggl)+
\sum_{I\subset I_0}(-1)^{|I|}\chi\biggl(\mathcal{D},\mathcal{G}(\underline{v})/\mathcal{G}(\underline{v}+\underline{1}_I)\biggr)=$$
$$=\sum_{I\subset
I_0}(-1)^{|I|}(h(\underline{v}-\underline{1}+\underline{1}_I)+h(\underline{v}))=
\sum_{I\subset
I_0}(-1)^{|I|}h(\underline{v}-\underline{1}+\underline{1}_I).$$
Otherwise, the coefficient at the monomial ${\underline{t}}^{\underline{v}}$
in the series $$\widetilde{P}'(\underline{t})=\biggl(\sum_{\underline{v}\in \mathbb{Z}^s}h(\underline{v})\cdot
\underline{t}^{\underline{v}}\biggr)\cdot\prod_{i=1}^s(t_i-1)$$ is also
equal to
$\sum_{I\subset
I_0}(-1)^{|I|}h(\underline{v}-\underline{1}+\underline{1}_I).$
\end{proof}
Denote $\mathcal{G}|_{E_i}$ by $\mathcal{G}_i$, and let $r$ be the rank of
$\mathcal{G}$.
\begin{lemma}
Suppose that
$0\rightarrow\mathcal{E}\rightarrow\mathcal{F}\rightarrow\mathcal{G}\rightarrow0$
is an exact sequence of sheaves on $\mathbb{P}^1$, $\mathcal{K}$ is a
locally free sheaf on $\mathbb{P}^{1}$. Then the sequence
$0\rightarrow\mathcal{K}\otimes\mathcal{E}\rightarrow\mathcal{K}\otimes\mathcal{F}\rightarrow\mathcal{K}\otimes\mathcal{G}\rightarrow0$
is also exact.
\end{lemma}
\begin{proof}
By the Birkhoff-Grothendieck splitting theorem every locally free sheaf on
 $\mathbb{P}^1$ is isomorphic to a direct sum of pervasive:
$\mathcal{K}=\sum_{i=1}^{r}\mathcal{K}_i$. Since the multiplication of the
exact triple on the pervasive sheaf does not disturb it exactness, sequences
$$0\rightarrow\mathcal{K}_i\otimes\mathcal{E}\rightarrow\mathcal{K}_i\otimes\mathcal{F}\rightarrow\mathcal{K}_i\otimes\mathcal{G}\rightarrow0$$
are exact. Therefore the desired sequence is exact as the direct sum of exact
ones.
\end{proof}
\begin{lemma}
Let $D=\sum_{i}k_{i}p_i$ be an effective divisor on
$\mathbb{P}^1$. If
$$0\rightarrow\mathcal{E}\rightarrow\mathcal{F}\rightarrow\mathcal{G}\rightarrow0$$
is an exact sequence of locally free sheaves, then
$$\chi(\mathcal{E}-D)-\chi(\mathcal{F}-D)+\chi(\mathcal{G}-D)=0.$$
\end{lemma}
\begin{proof}
$\chi(\mathcal{F}-D)=\chi(\mathcal{F})-deg(D)\cdot
rk(\mathcal{F})=\chi(\mathcal{E})+\chi(\mathcal{G})-deg(D)\cdot(rk(\mathcal{E})+rk(\mathcal{G}))=
\chi(\mathcal{E}-D)+\chi(\mathcal{G}-D).$
\end{proof}
\begin{lemma}
By $c_{(i)}$ denote the value of the first Chern class of $\mathcal{G}_i$ on
the fundamental homological class $[E_i]$. Then
$$\chi\biggl(E_i,\mathcal{G}(\underline{k})/\mathcal{G}(\underline{k}+\underline{1}_{\{i\}})|_{E_i}\biggr)=r+c_{(i)}+r\cdot\underline{k}\cdot\underline{d}^i.$$
\end{lemma}
\begin{proof}
Let $\nu_i$ be the normal sheaf to the $E_i$ in the variety $\mathcal{X}$.
One has an exact triple
$0\rightarrow\mathcal{J}_i^2\rightarrow\mathcal{J}_i\rightarrow\nu_{i}^{*}\rightarrow0$,
and since the normal bundle is one-dimensional,
the triple
$$0\rightarrow\mathcal{J}_i^{k_i+1}\rightarrow\mathcal{J}_i^{k_i}\rightarrow\nu_{i}^{-k_i}\rightarrow0$$
is also exact.
Since $\mathcal{G}_i$ is locally free, one has an exact sequence
$$0\rightarrow\mathcal{G}_i\otimes\mathcal{J}_i^{k_i+1}\rightarrow\mathcal{G}_i\otimes\mathcal{J}_i^{k_i}\rightarrow\mathcal{G}_i\otimes\nu_{i}^{-k_i}\rightarrow0,$$
therefore by lemma 4
$$\chi\biggl(E_i,\mathcal{G}(\underline{k}+\underline{1}_{\{i\}})/\mathcal{G}(\underline{k})\biggr)=\chi\biggl(E_i,\mathcal{G}_i\otimes\nu_{i}^{-k_i}\otimes\prod_{b\neq
i}\mathcal{J}_{b}^{k_b}\biggr).$$

Let $a$ be a canonical generator in  $H^2(\mathbb{P}^1,\mathbb{Z})$, then Chern character of $\mathcal{G}_i$
is equal to
$ch(\mathcal{G}_i)=r+c_{(i)}a$. Otherwise, $ch(\nu_{i}^{*})=1+d_{i}^{i}a$, then $$ch(\mathcal{G}_i\otimes\nu_{i}^{-k_i})=
(r+c_{(i)}a)(1+d_{i}^{i}a)^{k_i}=r+(c_{(i)}+rd_{i}^{i}k_i)a,$$ so $\chi(\mathcal{G}_i\otimes\nu_{i}^{-k_i})=r+c_{(i)}+rd_{i}^{i}k_i.$

Hence $$\chi\biggl(E_i,\mathcal{G}(\underline{k})/\mathcal{G}(\underline{k}+\underline{1}_{\{i\}})\biggr)=r+c_{(i)}+rd_{i}^{i}k_i-\sum_{E_j\cap
E_i\neq\emptyset}rk_j=r+c_{(i)}+r\sum_{j}d_{i}^{j}k_j.$$
\end{proof}
For $I\subset I_0$ by $\mu(I)$ denote a number of intersection points
of different lines with labels from the set $I$.
\begin{lemma}
Let $\mathcal{H}=\mathcal{G}(\underline{k})/\mathcal{G}(\underline{k}+\underline{1}_{I})$.
Then
$$\chi(\mathcal{D},\mathcal{H})=r\mu(I)+\sum_{i\in
I}\biggl(r+c_{(i)}+r(\underline{k}+\underline{1}_{I\setminus\{i\}})\underline{d}^i\biggr).$$
\end{lemma}
\begin{proof}
Consider $i\in I$. One could see that
$$\chi\biggl(E_i,\mathcal{G}(\underline{k})/\mathcal{G}(\underline{k}+\underline{1}_{I})\biggr)=$$
$$\chi\biggl(E_i,\mathcal{G}(\underline{k})/\mathcal{G}(\underline{k}+\underline{1}_{I\setminus\{i\}})\biggr)+
\chi\biggl(E_i,\mathcal{G}(\underline{k}+\underline{1}_{I\setminus\{i\}})/\mathcal{G}(\underline{k}+\underline{1}_{I})\biggr)=$$
$$-r\sum_{b\in
I\setminus\{i\}}d_{ib}+r+c_{(i)}+r(\underline{k}+\underline{1}_{I\setminus\{i\}})\underline{d}^i.$$
Since $\mathcal{H}$
is supported on $\cup_{i\in I}E_i$, from the Mayer-Vietoris exact sequence it follows that
 $$\chi(\mathcal{D},\mathcal{H})=\sum_{i\in I}\biggl(r+c_{(i)}+r(\underline{k}+\underline{1}_{I\setminus\{i\}})\underline{d}^i-r\sum_{b\in
I\setminus\{i\}}d_{ib}\biggr)-r\mu(I)=$$
$$=\sum_{i\in
I}\biggl(r+c_{(i)}+r(\underline{k}+\underline{1}_{I\setminus\{i\}})\underline{d}^i\biggr)+2r\mu(I)-r\mu(I).$$
\end{proof}
Using lemma 6, one could calculate the series $\widetilde{P}_{\mathcal{G}}(\underline{t})$
for
$\mathcal{G}=\mathcal{O}_{\mathcal{X}}$. Let $$m_i^{I}=-\underline{1}_{I\setminus \{i\}}\underline{d}^{i}$$ be the number of intersection points
$E_i$ with other components of exceptional divisor $\mathcal{D}$ with labels from the set $I$.
\begin{lemma}
If $s>2$ then $$\widetilde{P}_{\mathcal{O}}(\underline{t})\equiv 0.$$
\end{lemma}
\begin{proof}
From lemma 2 it follows that the coefficient at $\underline{t}^{\underline{k}}$ is
equal to
$$-\sum_{I\subset
I_0}(-1)^{|I|}\chi\biggl(\mathcal{D},\mathcal{G}(\underline{k})/\mathcal{G}(\underline{k}+\underline{1}_I)\biggr)=$$
$$-\sum_{I\subset I_0}(-1)^{|I|}\Biggl(\mu(I)+\sum_{i\in
I}\biggl(1+(\underline{k}+\underline{1}_{I\setminus\{i\}})\underline{d}^i\biggr)\Biggr).$$
Denote $v_i=1+\underline{k}\underline{d}^i,$ then the expression for the
coefficient could be represented in the form
 $$-\sum_{I\subset
I_0}(-1)^{|I|}(\mu(I)+\sum_{i\in I}(v_i-m_i^I)).$$ Since in the sum
$\sum_{i\in I}m_i^I$ every intersection point is taken into account twice,
it is clear that $\sum_{i\in I}m_i^I=2\mu(I)$. Hence the expression is equal
to $$-\sum_{I\subset I_0}(-1)^{|I|}\sum_{i\in I}(v_i-{m_i^I\over 2})=
-\sum_{i}(v_i\sum_{I\ni i}(-1)^{|I|}+\sum_{I\ni i}(-1)^{|I|}{m_i^I\over
2})=$$ $$=0+{1\over 2}\sum_{i}\sum_{I:i\in I,\forall j\in I
d_{i}^{j}<0}(-1)^{|I|}m_i^I\sum_{J:\forall j\in J d_{i}^{j}=0}(-1)^{|J|}. $$
If there are divisors, which does not intersect  $E_i$, one has
$$\sum_{J:\forall j\in J d_{i}^{j}=0}(-1)^{|J|}=0.$$ In the opposite case
$E_i$ is intersected by  $s-1$ divisors and the $i$-th term is equal to $-{1\over
2}\sum_{m=0}^{s-1}(-1)^{m}m{m\choose s-1}$ which vanishes if
$s-1>1.$
\end{proof}

This result is unnatural, but it shows that
 $P_{\mathcal{G}}$ could not be interpreted as generalization of the notion of Poincar\'e series in any sense.

Let $\psi(x)$ be equal to $x$ if $x\ge 0$ and equal to zero in the opposite
case.

\begin{lemma}
Let $\mathcal{G}=\mathcal{O}_{\mathcal{X}}.$ Then $$\dim
H^{0}\biggl(\mathcal{D},\mathcal{G}(\underline{k})/\mathcal{G}(\underline{k}+\underline{1}_I)\biggr)=
\mu(I)+\sum_{i\in
I}\psi\biggl(\chi(E_i,\mathcal{G}(\underline{k})/\mathcal{G}(\underline{k}+\underline{1}_I))\biggr).$$
\end{lemma}
\begin{proof}
Compute $\dim
H^{0}\biggl(E_i,\mathcal{G}(\underline{k}+\underline{1}_{I\setminus\{i\}})/\mathcal{G}(\underline{k}+\underline{1}_{I})\biggr).$
Let $f$ be the function on $\mathcal{X}$ representing the section of the corresponding sheaf on $E_i$.
Since
$\mathcal{J}_i^{k_i}/\mathcal{J}_i^{k_i+1}\simeq \nu_i^{-k_i}$, the class of $f$
is a well-defined section of
$\nu_i^{-k_i}\simeq\mathcal{O}(k_{i}d_{i}^{i})$. The space of global sections of this bundle which vanish at intersection points with other divisors
$E_j,j\neq i$ with order bigger or equal to
 $(\underline{k}+\underline{1}_{I\setminus\{i\}})_j$ is canonically isomorphic to the
$$H^{0}\biggl(\mathbb{CP}^1,\mathcal{O}(k_{i}d_{i}^{i}+\sum_{j\neq i}k_{j}d_{i}^{j}+\sum_{b\in I\setminus\{i\}}d_{i}^{j})\biggr)$$
and its dimension is equal to $\psi((\underline{k}+\underline{1}_{I\setminus\{i\}})\underline{d}^i+1).$
Furthermore, $$\dim
H^{0}\biggl(\mathbb{CP}^1,\mathcal{G}(\underline{k})/\mathcal{G}(\underline{k}+\underline{1}_{I\setminus\{i\}})\biggr)=
-\underline{1}_{I\setminus\{i\}}\underline{d}^i$$ and
$H^{1}\biggl(\mathbb{CP}^1,\mathcal{G}(\underline{k}+\underline{1}_{I\setminus\{i\}})/\mathcal{G}(\underline{k}+\underline{1}_{I})\biggr)=0$,
hence from the exact sequence
$$0\rightarrow\mathcal{G}(\underline{k}+\underline{1}_{I\setminus\{i\}})/\mathcal{G}(\underline{k}+\underline{1}_{I})\rightarrow
\mathcal{G}(\underline{k})/\mathcal{G}(\underline{k}+\underline{1}_{I})\rightarrow
\mathcal{G}(\underline{k})/\mathcal{G}(\underline{k}+\underline{1}_{I\setminus\{i\}})\rightarrow0$$
$$\dim H^{0}\biggl(E_i,\mathcal{G}(\underline{k})/\mathcal{G}(\underline{k}+\underline{1}_{I})\biggr)=
\psi\biggl(\chi(E_i,\mathcal{G}(\underline{k})/\mathcal{G}(\underline{k}+\underline{1}_{I}))\biggr)-\underline{1}_{I\setminus\{i\}}\underline{d}^i.$$

Since
$\mathcal{G}(\underline{k})/\mathcal{G}(\underline{k}+\underline{1}_I)$
is supported on the union of lines with numbers from $I$, it is sufficient to consider the space of its global sections over this union.
For such global section we could construct the collection of sections over
 $E_i$ which are uniquely determined by sets of their zeros up to the multiplication on a constant.

Otherwise, if we have a collection of points with multiplicities on
 $E_i$ , let us draw through every point a germ of analytical curve intersecting
 $E_i$ with the corresponding multiplicity. Under projection on
 $\mathbb{C}^2$ we will have a germ of a reducible curve.
Let us define this germ by equation
 $\{g=0\}$ and
consider a function $\pi^{*}g$ ($g$ is determined uniquely up to the multiplication on the
function,
which does not vanish at the origin).
Since on every $E_i$ section corresponding to $\pi^{*}g$
 has as much zeros as the section, from which we have started,
 it follows that $\pi^{*}g$
matches the same sections of the same powers of conormal bundles.

Therefore Mayer-Vietoris sequence is also exact in the part with global
sections. This note proves the lemma.
\end{proof}
Analogous to lemma 2 it is easy to prove that the coefficient at $\underline{t}^{\underline{k}}$
in $P_L$ is equal to $-\sum_{I\subset
I_0}(-1)^{|I|}\dim(L(\underline{k})/L(\underline{k}+\underline{1}_I)).$

\begin{definition}
The series
$$P^{g}_{\mathcal{G}}(\underline{t})=-\sum_{\underline{k}\in\mathbb{Z}^{r}_{\ge
0}}\underline{t}^{\underline{k}}\sum_{I\subset I_0}(-1)^{|I|}\biggl(
r\mu(I)+\sum_{i\in
I}\psi\left(\chi(E_i,\mathcal{G}(\underline{k}+\underline{1}_{I\setminus\{i\}})/\mathcal{G}(\underline{k}+\underline{1}_I))\right)\biggr).$$
is said to be geometrical Poincar\'e series of the filtration on the sheaf
$\mathcal{G}$.
\end{definition}

It is clear from lemma 7  that geometrical Poincar\'e series for
 $\mathcal{O}_{\mathbb{C}^2,0}$ coincides with the series
 $P_{H^{0}}$ computed in \cite{dg},
so it is reasonable to consider the geometrical Poincar\'e series as a
proper
generalization of the notion of Poicar\'e series for the space of global
sections.

\begin{lemma}
Let $a_i$ be nonnegative integers.
Then
$$-\sum_{I\subset I_0}(-1)^{|I|}\biggl(\mu(I)+\sum_{i:i\in
I,a_i+1-m_i^I>0}(a_i+1-m_{i}^{I})\biggr)=\chi(\prod_{i}S^{a_i}\widetilde{E}_i)).$$
\end{lemma}
\begin{proof}
Denote
$f(\underline{a},\mathcal{D})=\chi(\prod_{i}S^{a_i}\widetilde{E}_i).$ Let us
prove the proposition of lemma by induction on the number of exceptional lines.
If the line is unique, then
$a_1+1=(-1)^{a_1}{-2\choose a_1}=\chi(S^{a_1}E_1).$ If there are two lines, it's easy to check the lemma's proposition:
$(a_1+1)+(a_2+1)-((a_1+1-1)+(a_2+1-1)+1)=1.$ (if
$a_1>0, a_2>0;$ other cases are analogous).
Let us prove the inductive transition.
Since the dual graph of resolution is a tree, it has a vertex with degree 1,
i. e., a divisor $E_j$ which intersects the unique line $E_l$. Suppose that $a_l\neq 0.$ Тhen $\chi(\widetilde{E}_j)=1$,
hence $\chi(S^{a_j}\widetilde{E}_j)=1$; otherwise, if
$\hat{E_l}=\widetilde{E}_l\cup(E_j\cap E_l)$, then
$$S^{a_l}\hat{E_l}=S^{a_l}\widetilde{E}_l\sqcup(E_j\cap E_l)\times
S^{a_l-1}\hat{E_l},$$ and
$\chi(S^{a_l}\widetilde{E}_l)=\chi(S^{a_l}\hat{E_l})-\chi(S^{a_l-1}\hat{E_l}).$
Denote $\hat{\mathcal{D}}=\cup_{i\neq j}E_i, \widehat{m}_i^{I}=m_{i}^{I\setminus\{j\}}.$
Then
$\prod_{i}\chi(S^{a_i}\widetilde{E}_i)=\prod_{i\neq j}\chi(S^{a_i}\widetilde{E}_i)=f(\underline{a},\hat{\mathcal{D}})-f(\underline{a}-\underline{1}_{\{l\}},\hat{\mathcal{D}}).$
Otherwise, for $A,B,C\subset I_0$ denote $$\mathcal{M}_{A,B,C}=\{(i,I):i\in I\subset I_0,A\cap I=\emptyset,B\subset I,i\notin C,a_i+1-m_i^I>0\}.$$
$$\widehat{\mathcal{M}}_{A,B,C}=\{(i,I):i\in I\subset I_0,A\cap I=\emptyset,B\subset I,i\notin C,a_i+1-\widehat{m}_i^I>0\}.$$
Then
\begin{flushleft}
$\sum_{I\subset I_0}(-1)^{|I|}(\mu(I)+\sum_{i:i\in
I,a_i+1-m_i^I>0}(a_i+1-m_{i}^{I}))=\sum_{I}(-1)^{|I|}\mu(I)+\sum_{(i,I)\in\mathcal{M}_{\emptyset,\emptyset,\emptyset}}(-1)^{|I|}(a_i+1-m_i^I))=$
$\sum_{I:(j,I)\in\mathcal{M}_{\{l\},\emptyset,\emptyset}}(-1)^{|I|}(a_j+1-m_j^I)+ \sum_{I:(j,I)\in\mathcal{M}_{\emptyset,\{l\},\emptyset}}
(-1)^{|I|}(a_j+1-m_j^I)+ \sum_{I:(l,I)\in\mathcal{M}_{\emptyset,\{j\},\emptyset}}(-1)^{|I|}(a_l+1-m_l^I)+
\sum_{(i,I)\in\mathcal{M}_{\emptyset,\{j\},\{j,l\}}}(-1)^{|I|}(a_i+1-m_i^I)+ \sum_{I:j\in I,k\in
I}(-1)^{|I|}\mu(I)+\sum_{I:j\in I,k\notin I}(-1)^{|I|}\mu(I)-
f(\underline{a},\hat{\mathcal{D}})=\sum_{I:j\in I,l\notin
I}(-1)^{|I|}(a_j+1)+ \sum_{I:j\in I,l\in I,a_j>0}(-1)^{|I|}a_j+
\sum_{I:(l,i)\in\mathcal{M}_{\emptyset,\{j\},\emptyset}}(-1)^{|I|}((a_l-1)+1+\hat{m}_l^I)+$
$\sum_{(i,I)\in\widehat{\mathcal{M}}_{\emptyset,\{j\},\{j,l\}}}(-1)^{|I|}(a_i+1-\hat{m}_i^I)- \sum_{\hat{I}:j\notin
\hat{I},l\in \hat{I}}(-1)^{|\hat{I}|}(\mu(\hat{I})+1)-\sum_{\hat{I}:j\notin
\hat{I},l\notin \hat{I}}(-1)^{|\hat{I}|}\mu(\hat{I})-
f(\underline{a},\hat{\mathcal{D}})=0+0+f(\underline{a}-\underline{1}_{\{l\}},\hat{\mathcal{D}})-f(\underline{a},\hat{\mathcal{D}})=
-f(\underline{a},\mathcal{D}).$
\end{flushleft}
If $a_l=0$, then
$f(\underline{a},\mathcal{D})=f(\underline{a},\hat{\mathcal{D}}).$
Otherwise,
\begin{flushleft}
$\sum_{I}(-1)^{|I|}(\mu(I)+\sum_{i:i\in
I,a_i+1-m_i^I>0}(-1)^{|I|}(a_i+1-m_i^I))
=\sum_{I:(j,I)\in\mathcal{M}_{\{l\},\emptyset,\emptyset}}(-1)^{|I|}(a_j+1-m_j^I)+
\sum_{I:(j,I)\in\mathcal{M}_{\emptyset,\{l\},\emptyset}}(-1)^{|I|}(a_j+1-m_j^I)+
\sum_{(i,I)\in\mathcal{M}_{\emptyset,\{j\},\{j,l\}}}(-1)^{|I|}(a_i+1-m_i^I)+
\sum_{I:j\in I,l\in I}(-1)^{|I|}\mu(I)+\sum_{I:j\in I,l\notin
I}(-1)^{|I|}\mu(I)-
f(\underline{a},\hat{\mathcal{D}})=\sum_{I:j\in I,l\notin
I}(-1)^{|I|}(a_j+1)+
\sum_{I:j\notin I,k\in
I,a_j>0}(-1)^{|I|}a_j+
\sum_{(i,I)\in\widehat{\mathcal{M}}_{\emptyset,\{j,l\},\{j,l\}}}(-1)^{|I|}(a_i+1-\hat{m}_i^I)+$
$\sum_{(i,I)\in\widehat{\mathcal{M}}_{\{l\},\{j\},\{j,l\}}}(-1)^{|I|}(a_i+1-\hat{m}_i^I)-
\sum_{\hat{I}:j\notin \hat{I},l\in
\hat{I}}(-1)^{|\hat{I}|}\mu(\hat{I})-\sum_{\hat{I}:j\notin \hat{I},l\notin
\hat{I}}(-1)^{|\hat{I}|}\mu(\hat{I})-
f(\underline{a},\hat{\mathcal{D}})=0+0+$
$\sum_{(i,I\cup\{l\})\in\mathcal{M}_{\{j\},\{l\},\{j,l\}}}(-1)^{|I|}(a_i+1-\hat{m}_i^I+d_{il})-$

$\sum_{(i,I)\in\mathcal{M}_{\{j,\},\emptyset,\{j,l\}}}(-1)^{|I|}(a_i+1-\hat{m}_i^I)+
\sum_{(i,I\cup\{l\})\in\mathcal{M}_{\emptyset,\{j,l\},\{j,l\}}}(-1)^{|I|}(-d_{il})
-f(\underline{a},\hat{\mathcal{D}})=-f(\underline{a},\mathcal{D}).$
\end{flushleft}
\end{proof}
Let $\zeta_i=1$ if every exceptional line intersects $E_i$ and $\zeta_i=0$ in the opposite case.
\begin{lemma}
Suppose that $a_i$ are nonnegative integers and $u_i$ are arbitrary numbers.
Then
$$-\sum_{i,I:i\in
I,a_i+1-m_i^I>0}(-1)^{|I|}u_i=\sum_i\zeta_{i}u_{i}(-1)^{a_i}{1-\chi(\widetilde{E}_i)\choose
a_i}.$$
\end{lemma}
\begin{proof}
$-\sum_{i,I:i\in I,a_i+1-m_i^I>0}(-1)^{|I|}u_i=$
$$=-\sum_{i}u_i\sum_{I:i\in
I,a_i+1-m_i^I>0}(-1)^{|I|}
=$$
$$=\sum_{i}u_i\sum_{I:i\in I,\forall j\in I d_{i}^{j}<0,
a_i+1-|I|>0}(-1)^{|I|}\sum_{J:\forall j\in J d_{i}^{j}=0}(-1)^{|J|}.$$ The $i$-th term is equal to 0 if there exists
a line, which does not intersect $E_i$, and it is equal to $u_i\sum_{s=0}^{a_i}(-1)^{s}{2-\chi(\widetilde{E}_i)\choose s}=u_i(-1)^{a_i}{1-\chi(\widetilde{E}_i)\choose a_i},$
if all lines intersect $E_i$.
\end{proof}
\begin{theorem}
Geometrical Poincar\'e series of the filtration on the sheaf $\mathcal{G}$
is equal to the regular part of the Laurent series
$$r\prod_{i}{t_i}^{-\lceil{c_{(i)}\over
r}\rceil}(1-\underline{t}^{\underline{m}_i})^{-\chi(\widetilde{E}_i)}-
r\prod_{i}{t_i}^{-\lceil{c_{(i)}\over r}\rceil}(1-\underline{t}^{\underline{m}^{i}})^{-1}\sum_i\zeta_i\{-{с_i\over r}\}(1-\underline{t}^{\underline{m}^i})^{2-\chi(\tilde{E_i})},$$
where $\{x\}$ is a fractional part of $x$, and $\lceil{x}\rceil$ is the smallest integer bigger or equal to $x$.
\end{theorem}
\begin{proof}
From lemma 5 it follows that
the coefficient at $\underline{t}^{\underline{k}}$ in the geometrical Poincar\'e series
for sheaf $\mathcal{G}$ is equal to
$$-\sum_{I\subset I_0}(-1)^{|I|}\biggl(r\mu(I)+\sum_{i\in
I}\psi((r+c_{(i)}+r(\underline{k}+\underline{1}_{I\setminus\{i\}})\underline{d}^i))\biggr).$$
Denote $v_i=r+c_{(i)}+r\cdot\underline{k}\cdot\underline{d}^i$. It is clear that if there exists a term in the sum
which does not
vanish, then $v_i>0.$ Hence the coefficient is equal to
$$-\sum_{I\subset I_0}(-1)^{|I|}\biggl(r\mu(I)+\sum_{i\in
I,v_i-rm_i^{I}>0}(v_i-rm_i^{I})\biggr)=$$
$$-r\sum_{I\subset I_0}(-1)^{|I|}\biggl(\mu(I)+
\sum_{i\in I,{v_i\over r}-m_i^I>0}(\lceil{v_i\over
r}\rceil-\{-{v_i\over r}\}-m_i^I)\biggr)
=$$ $$-r\sum_{I\subset I_0}(-1)^{|I|}\biggl(\mu(I)+
\sum_{i\in I,\lceil{v_i\over r}\rceil-m_i^I>0}(\lceil{v_i\over
r}\rceil-m_i^I)\biggr)+
 r\sum_{i,I:i\in I,\lceil{v_i\over
r}\rceil-m_i^I>0}(-1)^{|I|}\{-{v_i\over r}\}=$$ $$r\chi\biggl(\prod_{i}S^{\lceil{v_i\over r}\rceil-1}\widetilde{E}_i\biggr)-
\sum_i\zeta_{i}\{-{v_i\over r}\}(-1)^{\lceil{v_i\over
r}\rceil-1}{1-\chi(\widetilde{E}_i)\choose \lceil{v_i\over r}\rceil-1}.$$

Furthermore, $\lceil{v_i\over r}\rceil=1+\lceil{c_{(i)}\over
r}\rceil+\underline{k}\underline{d}_i>0,$ i. e., $w_i=\lceil{c_{(i)}\over
r}\rceil+\underline{k}\underline{d}_i$ can be an arbitrary nonnegative
integer and
$$\underline{k}=(\underline{w}-\lceil{\underline{c}\over r}\rceil)M.$$
Therefore the Poincar\'e series is equal to
$$\sum_{\underline{w}\in\mathbb{Z}^{s}_{\ge
0}}\Biggl(r\chi\biggl(\prod_{i}S^{w_i}\widetilde{E}_i\biggr)-r\sum_i\zeta_{i}\{-{c_{(i)}\over r}\}(-1)^{w_i}{1-\chi(\widetilde{E}_i)\choose w_i}\Biggr)
\underline{t}^{(\underline{w}-\lceil{\underline{c}\over
r}\rceil)M}=$$
$$=r\prod_{i}{t_i}^{-\lceil{c_{(i)}\over
r}\rceil}\prod_{i}(1-t^{\underline{m}_i})^{-\chi(\widetilde{E}_i)}-$$ $$r\sum_i\zeta_{i}\{-{c_{(i)}\over r}\}\prod_{j}{t_j}^{-\lceil{c_j\over r}\rceil}
(1-\underline{t}^{\underline{m}^i})^{1-\chi(\widetilde{E}_i)}\prod_{j\neq i}(1-\underline{t}^{\underline{m}^j})^{-1}.$$
Generally speaking, this series contains also finite amount of terms with negative powers
$\underline{k}$.
We'll make the desired expression by throwing them away.
\end{proof}
\begin{theorem}
Geometrical Poincar\'e series for the sheaf of 1-forms is equal to the regular part of
the series
$$2\prod_{i}{t_i}^{-\lceil{d_{i}^{i}-2\over
2}\rceil}(1-\underline{t}^{\underline{m}_i})^{-\chi(\widetilde{E}_i)}-
2\prod_{i}{t_i}^{-\lceil{d_{i}^{i}-2\over 2}\rceil}(1-\underline{t}^{\underline{m}^{i}})^{-1}\sum_i\zeta_i\{-{d_{i}^{i}-2\over 2}\}(1-\underline{t}^{\underline{m}^i})^{2-\chi(\tilde{E_i})}.$$
\end{theorem}
\begin{proof}
In this case $\mathcal{G}=\Omega^{1}_{\mathcal{X}}.$ Then
$\mathcal{G}_i=T^{*}\mathcal{X}|_{E_i}.$ Consider an exact sequence
$$0\rightarrow TE_{i}\rightarrow
T\mathcal{X}|_{E_i}\rightarrow\nu_{i}\rightarrow0.$$
Dual complex is:
$$0\rightarrow\nu_i^{*}\rightarrow T^{*}\mathcal{X}|_{E_i}\rightarrow
T^{*}E_i\rightarrow 0.$$
Therefore $c_{(i)}=<с_1(T^{*}\mathcal{X}|_{E_i}),[E_i]>=d_{i}^{i}-2.$ Now
the proposition of the theorem follows from the theorem 1.
\end{proof}
If the line is unique, then series is equal to
${1+t\over (1-t)^2}.$ It differs from the Poincar\'e series on the space of global sections by 1.

If there are two divisors, the answer is equal to
${1+t_{1}t_2^2\over (1-t_{1}t_2)(1-t_{1}t_2^2)}.$

\vspace{1cm}

The author is grateful to S.M.Gusein-Zade for constant attention to
this work and useful discussions.

\end{document}